\documentclass[]{article}
\usepackage{amssymb}
\usepackage{amsmath}
\usepackage{amsthm}

\newtheorem{Theorem}{Theorem}
\newtheorem{Proposition}{Proposition}
\newtheorem{Definition}{Definition}

\newcommand{\pr}{\mbox{pr}_1\,}

\title {Stationary subsets of functional Menger $\cap$-algebras of multiplace functions}
\author {Wieslaw A. Dudek and Valentin S. Trokhimenko}
\date {}

\begin{document}
\sloppy \maketitle

\begin{abstract}
A functional Menger $\cap$-algebra is a set of $n$-place functions
containing $n$ projections and closed under the so-called Menger's
compositions of $n$-place functions and the set-theoretic
intersection of functions. We give the abstract characterization
for these subsets of functional Menger $\cap$-al\-gebras which
contain functions with fixed points.
\end{abstract}

\section{Introduction}
Investigation of multiplace functions by algebraic methods plays a
very important role in modern mathematics where we consider
various operations on sets of functions which are naturally
defined. The basic operation for functions is superposition
(composition), but there are  some other naturally defined
operations, which are also worth of consideration. For example,
the operation of set-theoretic intersection and the operation of
projections  (see for example \cite{Dudtro, Dudtro1, Dudtro2,
schtr, Tr-1}). In this study the central role play sets of
functions with fixed points. The study of such sets for functions
of one variable was initiated by B. M. Schein in \cite{146} and
\cite{SchLectures}. Then it was continued by V. S. Trokhimenko
(see \cite{Tr-2, Tr-3, Tr-4}) and W.A. Dudek \cite{Dudtro3} for
sets of functions of $n$ variables.

In this paper, we consider the sets of $n$-place functions
containing $n$-projections and closed under the so-called Menger's
composition and set-theoretic intersection of $n$-place functions.
For such functional Menger $\cap$-al\-gebras we give the abstract
characterization for subsets of functions with fixed points.

\section{Preliminaries}

Let $A^n$ be the $n$-th Cartesian product of a set $A$. Any
partial mapping from $A^n$ into $A$ is called an \textit{$n$-place
function} on $A$. The set of all such mappings is denoted by
$\mathcal{F}(A^n,A)$. On ${\mathcal F}(A^n,A)$ we define one
$(n+1)$-ary superposition $\mathcal{O}\colon (f,g_1,\ldots,
g_n)\mapsto f[g_1\ldots g_n]$, called the {\it Menger's
composition}, and $n$ unary operations $\mathcal{R}_i\colon
f\mapsto\mathcal{R}_if$, $i\in\overline{1,n}=\{1,\ldots,n\}$
putting
\begin{eqnarray}
f [g_1\ldots g_n](a_1,\ldots,a_n)=f(g_1(a_1,\ldots,a_n),\ldots,
g_n(a_1,\ldots,a_n)), \label{super} \\[4pt]
\mathcal{R}_if(a_1,\ldots,a_n)=a_i, \ \mbox{where} \ \
(a_1,\ldots,a_n) \in\pr f,
\end{eqnarray}
for $f,g_1,\ldots,g_n\in\mathcal{F}(A^n,A)$, $(a_1,\ldots,a_n)\in
A^n$, where $\pr f$ denotes the domain of a function $f$. It is
assumed the left and right hand side of equality (\ref{super}) are
defined, or not defined, simultaneously. Algebras of the form
$(\Phi,\mathcal{O},\cap,\mathcal{R}_1,\ldots,\mathcal{R}_n)$,
where $\Phi\subseteq\mathcal{F}(A^n, A)$ and $\cap$ is a
set-theoretic intersection, are called \textit{functional Menger
$\cap$-algebras of $n$-place functions}. In the literature such
algebras are also called \textit{functional Menger
$\cal{P}$-algebras} (see \cite{Dudtro} and \cite{Tr-4}). The set
of functions from $\Phi$ for which there exists fixed point, i.e.,
the set
 \[
\mathbf{St}(\Phi)= \{f\in\Phi\,|\,(\exists a\in
A)\,f(a,\ldots,a)=a\},
 \]
is called the \textit{stationary subset} of $\Phi$.

Let $(G,o)$ be a nonempty set with one $(n+1)$-ary operation
$$o\colon\ (x_0,x_1,\ldots,x_n)\mapsto x_0[x_1\ldots x_n].$$
An algebra $\mathcal{G}=(G,o,\curlywedge, R_1,\ldots,R_n)$ of type
$(n+1,2,1,\ldots,1)$ where $(G,\curlywedge)$ is a semilattice and
for all $i,k\in\overline{1,n}$ it satisfies the following axioms:
\[\begin{array}{cl}
\mathbf{A_1}\colon& \quad x[y_1\ldots y_n][z_1\ldots z_n]=
x[y_1[z_1\ldots z_n]\ldots y_n[z_1\ldots z_n]], \\[4pt]
\mathbf{A_2}\colon& \quad x[R_1x\ldots R_nx]=x, \\[4pt]
\mathbf{A_3}\colon& \quad x[\bar{u}\,|_iz][R_1y\ldots R_ny]=
x[\bar{u}\,|_iz[R_1y\ldots R_ny]], \\[4pt]
\mathbf{A_4}\colon& \quad R_ix[R_1y\ldots R_ny]= (R_ix)[R_1y\ldots
R_ny], \\[4pt]
\mathbf{A_5}\colon& \quad x[R_1y\ldots R_ny][R_1z\ldots R_nz]=
x[R_1z\ldots R_nz][R_1y\ldots R_ny], \\[4pt]
\mathbf{A_6}\colon& \quad R_ix[y_1\ldots y_n]=R_i((R_kx)[y_1\ldots y_n]), \\[4pt]
\mathbf{A_7}\colon& \quad (R_ix)[y_1\ldots y_n]=
y_i[R_1(x[y_1\ldots y_n])\ldots R_n(x[y_1\ldots y_n])],\\[4pt]
\mathbf{A_8}\colon& \quad x\curlywedge y[R_1z\ldots R_nz]=
(x\curlywedge y)[R_1z\ldots R_nz], \\[4pt]
\mathbf{A_9}\colon& \quad x\curlywedge y=x[R_1(x\curlywedge
y)\ldots R_n(x\curlywedge y)], \\[4pt]
\mathbf{A_{10}}\colon& \quad (x\curlywedge y)[z_1\ldots z_n]=
x[z_1\ldots z_n]\curlywedge y[z_1\ldots z_n],
\end{array}
\]
where $x[\bar{u}\,|_iz]$ means $x[u_1\ldots
u_{i-1}z\,u_{i+1}\ldots u_n]$, is called a \textit{functional
Menger $\curlywedge$-algebra of rank $n$}.

Any Menger algebra of rank $n$, i.e., an abstract groupoid $(G,o)$
with an $(n+1)$-ary operation satisfying $\mathbf{A_1}$, is
isomorphic to some set of $n$-place functions closed under
Menger's composition \cite{schtr}. Functional Menger
$\curlywedge$-algebras  are isomorphic to some functional Menger
$\cap$-algebras  of $n$-place functions (see \cite{Dudtro}). Each
homomorphism of such abstract algebras into corresponding algebras
of $n$-place functions is called a \textit{representation by
$n$-place functions}. Representations which are isomorphisms are
called \textit{faithful}.

Let $(P_i)_{i\in I}$ be the family of representations of a Menger
algebra $(G,o)$ of rank $n$ by $n$-place functions defined on sets
$(A_i)_{i\in I}$, respectively. By the \textit{union} of this
family we mean the mapping $P \colon \, g \mapsto P(g)$, where
$g\in G$, and $P(g)$ is an $n$-place function on $A= \bigcup
\limits_{i\in I}A_i$ defined by
 \[
P(g)=\bigcup\limits_{i\in I}P_i(g).
 \]
If $A_i\cap A_j=\emptyset$ for all $i,j\in I$, $i\neq j$, then $P$
is called the \textit{sum} of $(P_i)_{i\in I}$ and is denoted by
$P=\sum\limits_{i\in I}P_i$. It is not difficult to see that the
sum of representations is a representation, but the union of
representations may not be a representation (see for example
\cite{Dudtro} -- \cite{Tr-1}).

Let $\mathcal{G}=(G,o,\curlywedge, R_1,\ldots,R_n)$ be an
functional Menger $\curlywedge$-algebra of rank $n$. We shall say
 that a nonempty subset $H$ of $G$ is called
\begin{itemize}
\item \textit{$\curlywedge$-quasi-stable}, if for all $x\in G$
\[
x\in H\longrightarrow x[x\ldots x]\curlywedge x\in H,
\]
\item an \textit{$l$-ideal}, if for all $x,y_1,\ldots,y_n\in
G$
\[
(y_1,\ldots,y_n)\in G^{\,n}\setminus (G\setminus
H)^n\longrightarrow x[y_1\ldots y_n]\in H.
\]
\end{itemize}
A binary relation $\rho\subseteq G\times G$ is called
\begin{itemize}
\item \textit{stable}, if
\[
(x,y),(x_1,y_1),\ldots,(x_n,y_n)\in\rho\longrightarrow
(x[x_1\ldots x_n],y [y_1\ldots y_n])\in\rho
\]
for all $x,y,x_i,y_i\in G$, $i\in\overline{1,n}$,
\item \textit{$l$-regular}, if
\[
(x,y)\in\rho\longrightarrow (x[z_1\ldots z_n],y[z_1\ldots
z_n])\in\rho
\]
for all $x,y,z_i\in G$, $i\in\overline{1,n}$,
\item \textit{$v$-regular}, if
\[
(x_1,y_1),\ldots, (x_n,y_n)\in\rho\longrightarrow (z[x_1\ldots
x_n],z[y_1\ldots y_n])\in\rho
\]
for all $x_i,y_i,z\in G$, $i\in\overline{1,n}$,
\item \textit{$i$-regular}, where $i\in\overline{1,n}$, \ if
\[
(x,y)\in\rho\longrightarrow (u[\bar{w}\,|_ix],
u[\bar{w}\,|_iy])\in\rho
\]
for all $x,y,u\in G$, $\bar{w}\in G^n$,
\item \textit{$v$-negative}, if
\[
 (x[y_1\ldots y_n],y_i)\in\rho
\]
for all $x,y_1,\ldots,y_n\in G$ and $i\in\overline{1,n}$.
\end{itemize}
On $\mathcal{G}$ we define two binary relations $\zeta$ and $\chi$
putting
\[
(x,y)\in\zeta\longleftrightarrow x=y[R_1x\ldots R_nx], \qquad
(x,y)\in\chi \longleftrightarrow (R_1x,R_1y)\in\zeta .
\]
The first relation is a stable order, the second is an $l$-regular
and $v$-negative quasi-order containing $\zeta$ (see \cite{Tr-1}).
For these two relations the following conditions are valid:
\[
\begin{array}{lll}
x\leqslant y\longrightarrow R_ix\leqslant R_iy, \ \ \
i\in\overline{1,n}, &&x\sqsubset y\longleftrightarrow
R_ix\leqslant R_iy, \ \ \ i\in\overline{1,n},
\\[4pt]
x\sqsubset y\longleftrightarrow x[R_1y\ldots R_ny]=x,&&
(R_ix)[y_1\ldots y_n]\leqslant y_i, \ \ \ i\in\overline{1,n},
\\[4pt]
x[R_1y_1\ldots R_ny_n]\leqslant x, &&R_ix=R_iR_kx, \ \ \
i,k\in\overline{1,n},
\end{array}
\]
where $x\leqslant y\longleftrightarrow (x,y)\in\zeta$, and
$x\sqsubset y\longleftrightarrow (x,y)\in\chi$.

Let $W$ be the empty set or an $l$-ideal which is an
$\mathcal{E}$-class of a $v$-regular equivalence relation
$\mathcal{E}$ defined on a Menger algebra $(G,o)$ of rank $n$.
Denote by $(H_a)_{a\in A_{\mathcal{E}}}$ the family of all
$\mathcal{E}$-classes (uniquely indexed by elements of some set
$A_{\mathcal{E}}$) such that $H_a\ne W$. Next, for every $g\in G$
we define on $A_{\mathcal{E}}$ an $n$-place function
$P_{(\mathcal{E},W)}(g)$ putting
\begin{equation} \label{simplrep}
P_{(\mathcal{E},W)}(g)(a_1,\ldots,a_n)=b\longleftrightarrow
g[H_{a_1}\ldots H_{a_n}]\subseteq H_b,
\end{equation}
where $(a_1,\ldots,a_n)\in\pr
P_{(\mathcal{E},W)}(g)\longleftrightarrow g[H_{a_1}\ldots
H_{a_n}]\cap W=\emptyset$, and $H_b$ is an $\mathcal{E}$-class
containing all elements of the form $g[h_1\ldots h_n]$,\, $h_i\in
H_{a_i}$, $i\in\overline{1,n}$. It is not difficult to see
\cite{Dudtro2} that the mapping
$P_{(\mathcal{E},W)}\colon\,g\mapsto P_{(\mathcal{E},W)}(g)$
satisfies the identity
\begin{equation}\label{srep}
  P_{(\mathcal{E},W)}(g[g_1\ldots g_n])=
P_{(\mathcal{E},W)}(g)[P_{(\mathcal{E},W)}(g_1)\ldots
P_{(\mathcal{E},W)}(g_n)],
\end{equation}
i.e., $P_{(\mathcal{E},W)}$ is a representation of $(G,o)$ by
$n$-place functions. This representation will be called
\textit{simplest}.

\section{Stationary subsets}

The important properties of the stationary subset
$\mathbf{St}(\Phi)$ of the algebra $(\Phi,\mathcal{O})$, where
$\Phi\subseteq\mathcal{F}(A^n, A)$, are presented in the following
proposition.

\begin{Proposition} \label{P-1}
The stationary subset $\mathbf{St}(\Phi)$ of the algebra
$(\Phi,\mathcal{O})$ has the following properties:
\begin{eqnarray}
\label{st-1} f\subseteq g\wedge f\in\mathbf{St}(\Phi)\longrightarrow g\in\mathbf{St}(\Phi),\\[4pt]
\label{st-2} f\in\mathbf{St}(\Phi)\longrightarrow f[f\ldots f]\in\mathbf{St}(\Phi),\\[4pt]
\label{st-3} f\in\mathbf{St}(\Phi)\longrightarrow\mathcal{R}_if\in\mathbf{St}(\Phi),\\[4pt]
\label{st-4} f[g\ldots g]=g\wedge g\in\mathbf{St}(\Phi)\longrightarrow f\in\mathbf{St}(\Phi),\\[4pt]
\label{st-5} f[g\ldots g]=g\neq\mathbf{0}\longrightarrow f\in\mathbf{St}(\Phi),\\[4pt]
\label{st-6} f[g\ldots g]\neq\mathbf{0}\longrightarrow\mathcal{R}_if\in\mathbf{St}(\Phi),\\[4pt]
\label{st-7} f[g\ldots g]\cap g\neq\mathbf{0}\longrightarrow f\in\mathbf{St}(\Phi), \\[4pt]
\label{st-8}
\mathbf{0}\not\in\mathbf{St}(\Phi)\longrightarrow\mathcal{R}_i\mathbf{0}=\mathbf{0},
\end{eqnarray}
where $i\in\overline{1,n}$, $f,g\in\Phi$ and $\mathbf{0}$ is a
zero of $(\Phi,\mathcal{O})$.
\end{Proposition}
\begin{proof} If $f\in\mathbf{St}(\Phi)$, then $f(a,\ldots,a)=a$, whence, by $f\subseteq g$, we obtain
$g(a,\ldots,a)=a$. Thus $g\in\mathbf{St}(\Phi)$. This proves
(\ref{st-1}).

For $f\in\mathbf{St}(\Phi)$ we have also
\[
f[f\ldots f](a,\ldots,a)=f(f(a,\ldots,a),\ldots,f(a,\ldots,a))=
f(a,\ldots,a)=a.
\]
This implies $f[f\ldots f]\in\mathbf{St}(\Phi)$. So, (\ref{st-2})
is valid too. The proof of (\ref{st-3}) is analogous.

If $g\in\mathbf{St}(\Phi)$ and $f[g\ldots g]=g$, then, for some
$a\in A$, we have $g(a,\ldots,a)=a$. 
Consequently, $f(a,\ldots,a)=f(g(a,\ldots,a),\ldots,
g(a,\ldots,a))=f[g\ldots g](a,\ldots,a)=g(a,\ldots,a)=a$. Thus
$f\in\mathbf{St}(\Phi)$. This proves (\ref{st-4}).

Let now $f[g\ldots g]=g\neq\mathbf{0}$, where $\mathbf{0}$ is a
zero of $(\Phi,\mathcal{O})$. Then $g\neq\emptyset$. Thus there
exists $\bar{a}=(a_1,\ldots,a_n)\in\pr g$. Therefore
$f(g(\bar{a})\ldots g(\bar{a}))=f[g\ldots g](\bar{a})=g(\bar{a})$,
which implies $f\in\mathbf{St}(\Phi)$. The condition (\ref{st-5})
is proved. Similarly we can prove (\ref{st-6}) and (\ref{st-7}).

Observe that
\begin{equation} \label{null}
 \mathbf{0}\neq\emptyset\longleftrightarrow\mathbf{St}(\Phi)=\Phi.
\end{equation}
Indeed, if $\mathbf{0}\neq\emptyset$, then $\mathbf{0}(\bar{a})=b$
for some $\bar{a}\in A^n$ and $b\in A$. Since
$\mathbf{0}=f[\mathbf{0}\ldots\mathbf{0}]$ for every $f\in\Phi$,
we have
$b=\mathbf{0}(\bar{a})=f(\mathbf{0}(\bar{a}),\ldots,\mathbf{0}(\bar{a}))=
f(b,\ldots, b)$. Thus $f\in\mathbf{St}(\Phi)$. Consequently,
$\mathbf{St}(\Phi)=\Phi$. Conversely, if $\mathbf{St}(\Phi)=\Phi$,
then $\mathbf{0}\in\mathbf{St}(\Phi)$. Therefore
$\mathbf{0}(a,\ldots,a)=a$ for some $a\in A$. So,
$\mathbf{0}\neq\emptyset$.

Using just proved equivalence we can see that in the case
$\mathbf{0}\not\in\mathbf{St}(\Phi)$, i.e.,
$\mathbf{St}(\Phi)\neq\Phi$, must be $\mathbf{0}=\emptyset$.
Therefore
$\mathcal{R}_i\mathbf{0}=\mathcal{R}_i\emptyset=\emptyset
=\mathbf{0}$ for every $i\in\overline{1,n}$. This proves
(\ref{st-8}) and completes the proof of Proposition~\ref{P-1}.
\end{proof}
Note that in a functional Menger $\cap$-algebra
$(\Phi,\mathcal{O},\cap,\mathcal{R}_1,\ldots,\mathcal{R}_n)$ of
$n$-place functions without a zero $\mathbf{0}$ the subset
$\mathbf{St}(\Phi)$ coincides with $\Phi$. Indeed, in this algebra
$f\neq\emptyset$ for any $f\in\Phi$. Consequently, $f\cap
g[f\ldots f]\neq\emptyset$ for all $f,g\in\Phi$. Hence, $g[f\ldots
f](\bar{a})=f(\bar{a})$, i.e.,
$g(f(\bar{a}),\ldots,f(\bar{a}))=f(\bar{a})$ for some $\bar{a}\in
A^n$. This means that $g\in\mathbf{St}(\Phi)$. Thus,
$\Phi\subseteq\mathbf{St}(\Phi)\subseteq\Phi$, i.e.,
$\mathbf{St}(\Phi)=\Phi$. Therefore, below we will consider only
functional Menger $\cap$-algebras with a zero.

Let $\mathcal{G}=(G,o,\curlywedge, R_1,\ldots,R_n)$ be a
functional Menger $\curlywedge$-algebra of rank $n$.
\begin{Definition} \label{D-2}\rm
A nonempty subset $H$ of $G$ is called a \textit{stationary
subset} of a functional Menger $\curlywedge$-algebra $\mathcal{G}$
of rank $n$ if there exists its faithful representation $P$ by
$n$-place functions such that
\begin{equation} \label{stac}
 g\in H\longleftrightarrow P(g)\in\mathbf{St}(P(G))
\end{equation}
for every $g\in G$, where $P(G)=\{P(g)\,|\,g\in G\}$.
\end{Definition}

\begin{Theorem} \label{T-5a}
For a nonempty subset $H$ of $G$ to be a stationary subset of a
functional Menger $\curlywedge$-algebra $\mathcal{G}$ with a zero
$\mathbf{0}$, it is necessary and sufficient to be a
$\curlywedge$-quasi-stable subset satisfying for all $x,y\in G$
and $i\in\overline{1,n}$ the following three conditions:
\begin{eqnarray}
 \label{T5-1}\mathbf{0}\not\in H\longrightarrow R_i\mathbf{0}=\mathbf{0}, \\[4pt]
 \label{T5-2} x[y\ldots y]=y\in H\longrightarrow x\in H, \\[4pt]
 \label{T5-3} x[y\ldots y]\curlywedge y\neq\mathbf{0}\longrightarrow x\in
H.
\end{eqnarray}
\end{Theorem}
\begin{proof} The necessity of these conditions is a consequence
of our Proposition~\ref{P-1}, therefore we shall prove only their
sufficiency. For this assume that a nonempty subset $H$ of $G$
satisfies all the conditions of the theorem and consider on $G$ a
binary relation $\mathcal{E}_g$ and a subset $W_g$ defined in the
following way
\[\begin{array}{c}
\mathcal{E}_g=\{(x,y)\in G\times G\,|\,g\sqsubset x\curlywedge y
\vee
x,y\not\in\chi\langle g\rangle\}, \\[4pt]
W_g=\{x\in G\,|\,x\not\in\chi\langle g\rangle\},
\end{array}
\]
where $g\in G$ is fixed. Such defined relation $\mathcal{E}_g$ is
a $v$-regular equivalence for which each nonempty $W_g$ is an
$\mathcal{E}_g$-class and an $l$-ideal simultaneously (for details
see \cite{Dudtro}). Thus the pair $(\mathcal{E}_g,W_g)$ determines
the simplest representation $P_g=P_{(\mathcal{E}_g,W_g)}$ of
$(G,o)$ by $n$-place functions. In \cite{Dudtro} it is proved,
that $P_g$ is also a representation of $\mathcal{G}$, because it
satisfies, except (\ref{srep}), the equalities
\[
P_g(x\curlywedge
y)=P_g(x)\cap P_g(y),\quad P_g(R_ix)=\mathcal{R}_iP_g(x)
\]
for all $x,y\in G$ and $i\in\overline{1,n}$. Hence
$P=\sum\limits_{h\in G_0}P_h$, where
$G_0=G\setminus\{\mathbf{0}\}$ if $\mathbf{0}\not\in H$, and
$G_0=G$, if $\mathbf{0}\in H$, is a  representation of
$\mathcal{G}$ also. We must only prove that $P$ is faithful
representation, which satisfies the condition (\ref{stac}).

First we shall show that $H$ satisfies the condition
\begin{equation}
  \label{st-2a} \mathbf{0}\in H\longrightarrow H=G.
\end{equation}
Indeed, let $g\in G$ be any element of algebra $\mathcal{G}$, then
$g[\mathbf{0}\ldots \mathbf{0}]=\mathbf{0}\in H$, from where by
axiom (\ref{T5-2}) we obtain $g\in H$. So, $G\subseteq H\subseteq
G$, hence, $H=G$. The condition (\ref{st-2a}) is proved.

Now we shall prove that $H$ is its stationary subset of
$\mathcal{G}$. Let $g\in H$ and $P_0=P_{(\mathcal{E}_0,W_0)}$. If
$\mathbf{0}\in H$, then $G_0=G=H$, whence $\mathbf{0}\in G_0$.
Since $g\not\in W_g$ for every $g\in G$, we have
$\mathbf{0}\not\in W_{\mathbf{0}}$. Let $X$ be this
$\mathcal{E}_{\mathbf{0}}$-class, indexed by $a$, which contains
$\mathbf{0}$. Clearly $X\neq W_{\mathbf{0}}$. From
$g[\mathbf{0}\ldots\mathbf{0}]=\mathbf{0}$, applying the
$v$-regularity of $\mathcal{E}_{\mathbf{0}}$, we obtain $g[X\ldots
X]\subseteq X$. Consequently, $P_{\mathbf{0}}(g)(a,\ldots,a)=a$.
Hence, $P(g)(a,\ldots,a)=a$, which proves
$P(g)\in\mathbf{St}(P(G))$.

Now let $\mathbf{0}\not\in H$. Then $G_0=G\setminus\{\mathbf{0}\}$
and $h=g\curlywedge g[g\ldots g]\in H$ for every $g\in H$, because
$H$ is $\curlywedge$-quasi-stable. Thus $h\neq\mathbf{0}$, whence
$h\in G_0$. We shall consider the representation $P_h$. Since
$h\not\in W_h$, we have $g\curlywedge g[g\ldots g]\not\in W_h$.
Consequently, $g\equiv g[g\ldots g](\mathcal{E}_h)$. Moreover,
from
\[ \arraycolsep=.5mm\begin{array}{ll}
 h[R_1g\ldots R_ng]&=(g\curlywedge g[g\ldots g])[R_1g\ldots
R_ng]\\[4pt]
 &\stackrel{\mathbf{A}_8}{=}g\curlywedge g [g\ldots g] [R_1g\ldots R_ng] \\[4pt]
&\stackrel{\mathbf{A}_1}{=}g\curlywedge g[g[R_1g\ldots R_ng]\ldots g [R_1g\ldots R_ng]] \\[4pt]
&\stackrel{\mathbf{A}_2}{=}g\curlywedge g[g\ldots g]=h
\end{array}
\]
it follows $h [R_1g\ldots R_ng]=h$. Therefore $h\sqsubset g$,
which means that $g\not\in W_h$. Let $Y$ denotes the
$\mathcal{E}_h$-class containing $g$. Clearly, $Y\neq W_h$ and
$g[Y\ldots Y]\subseteq Y$. Hence $P_h(g)(b,\ldots,b)=b$, where $b$
is an element used as index of $Y$. Thus $P(g)(b,\ldots,b)=b$.
This means that also in this case $g\in H$ implies
$P(g)\in\mathbf{St}(P(G))$.

To prove the converse implication let $P(g)\in\mathbf{St}(P(G))$
for some $g\in G$. Because $P=\sum\limits_{h\in G_0}P_h$, there
exists $h\in G$ such that $P_{h}(g)$ has a fixed point. If
$\mathbf{0}\not\in H$, then $G_0=G\setminus\{\mathbf{0}\}$ and
$h\neq \mathbf{0}$. Let $X=H_a$ be this $\mathcal{E}_h$-class for
which $P_h(g)(a,\ldots,a)=a$, i.e., $g[X\ldots X]\subseteq X$,
where $X\neq W_h$. Obviously, for any $x\in X$ we have $g[x\ldots
x]\equiv x(\mathcal{E}_h)$. This means that $x\curlywedge
g[x\ldots x]\not\in W_h$ for any $x\in X$. Therefore
\[
h[R_1(x\curlywedge g[x\ldots x])\ldots R_n(x\curlywedge g[x\ldots
x])]=h\neq\mathbf{0},
\]
whence $R_i(x\curlywedge g[x\ldots x])\neq\mathbf{0}$ for every
$i\in\overline{1,n}$. This, in view of (\ref{T5-1}), gives
$x\curlywedge g[x\ldots x]\neq\mathbf{0}$. In the opposite case we
have $\mathbf{0}\in H$, which is impossible. Applying (\ref{T5-3})
to $x\curlywedge g[x\ldots x]\neq\mathbf{0}$ we obtain $g\in H$.
If $\mathbf{0}\in H$, then by the condition (\ref{st-2a}) we have
$H=G$ and therefore $g\in H$.

Thus we have proved that $g$ satisfies (\ref{stac}). So, $H$ is a
stationary subset of $\mathcal{G}$.

For completeness of the proof we must show that the representation
$P$ is faithful. If $P(g_1)=P(g_2)$, then $P(g_1)\subseteq P(g_2)$
and $P(g_2)\subseteq P(g_1)$, whence $g_1\leqslant g_2$ and
$g_2\leqslant g_1$ (for details see \cite{Dudtro}). This implies
$g_1=g_2$, because $\leqslant$ is an order.
 \end{proof}

Conditions formulated in the above theorem are not identical with
the conditions used for a characterization of stationary subsets
of restrictive Menger $\mathcal P$-algebras (see Theorem 8 in
\cite{Tr-4}). For example, the implication $$x\leqslant y\wedge
x\in H\longrightarrow y\in H$$ is omitted. Nevertheless, as it is
proved below, stationary subsets of functional Menger
$\curlywedge$-algebras with a zero have the same properties as
stationary subsets of restrictive Menger $\mathcal P$-algebras.

\begin{Theorem}\label{T-prop}
For a stationary subset $H$ of a functional Menger
$\curlywedge$-algebra $\mathcal{G}$ with a zero $\mathbf{0}$ the
following implications:
\begin{eqnarray}
\label{st-1a} \mathbf{0}\not\in H\longrightarrow \mathbf{0}\leqslant x,&&  \\[4pt]
\label{53}  x\leqslant y\wedge x\in H\longrightarrow y\in H,&&  \\[4pt]
\label{56} x\in H\longrightarrow x[x\ldots x]\in H, && \\[4pt]
\label{57} x\in H\longrightarrow R_ix\in H, && \\[4pt]
\label{59} x[y\ldots y]\neq\mathbf{0}\longrightarrow R_ix\in H, && \\[4pt]
\label{58} x[y\ldots y]=y\neq\mathbf{0}\longrightarrow x\in H, && \\[4pt]
\label{60} x\in H\wedge x\sqsubset y\longrightarrow R_iy\in H, && \\[4pt]
\label{61} \mathbf{0}\not\in H\wedge
x\sqsubset\mathbf{0}\longrightarrow x=\mathbf{0}
\end{eqnarray}
are valid for all $x,y\in G$ and $i=\overline{1,n}$.
\end{Theorem}
\begin{proof}
If $\mathbf{0}\not\in H$, then, by (\ref{T5-1}), we obtain
$R_i\mathbf{0}=\mathbf{0}$ for all $i\in\overline{1,n}$. Hence,
$\mathbf{0}=x[\mathbf{0}\ldots\mathbf{0}]=x[R_1\mathbf{0}\ldots
R_n\mathbf{0}]$, i.e., $\mathbf{0}=x[R_1\mathbf{0}\ldots
R_n\mathbf{0}]$ for any $x\in G$. So, $\mathbf{0}\leqslant x$.
This proves (\ref{st-1a}).

Now, let the premise of (\ref{53}) be satisfied, i.e., $x\leqslant
y$ and $x\in H$ for some $x,y\in G$. If $\mathbf{0}\in H$, then
$H=G$ by (\ref{st-2a}). Therefore $y\in H$. If $\mathbf{0}\not\in
H$, then, according to (\ref{st-1a}), for every $x\in G$ we have
$\mathbf{0}\leqslant x$. Since $x\in H$, the
$\curlywedge$-quasi-stability of $H$ implies $x[x\ldots
x]\curlywedge x\in H$. Hence $x[x\ldots x]\curlywedge
x\neq\mathbf{0}$, because $\mathbf{0}\not\in H.$ From $x\leqslant
y$, by the stability of $\leqslant$, we conclude $x[x\ldots
x]\leqslant y[x\ldots x]$. Consequently, $x[x\ldots x]\curlywedge
x\leqslant y[x\ldots x]\curlywedge x$. However
$\mathbf{0}\leqslant x[x\ldots x]\curlywedge x$, therefore
$$
\mathbf{0}\leqslant x[x\ldots x]\curlywedge x\leqslant y[x\ldots
x]\curlywedge x.
$$
Since $x[x\ldots x]\curlywedge x\neq\mathbf{0}$, the above gives
$y[x\ldots x]\curlywedge x\neq\mathbf{0}$, because in the opposite
case, by antisymmetry of $\leqslant$, we obtain $x[x\ldots
x]\curlywedge x=\mathbf{0}$, which is impossible. So, $y[x\ldots
x]\curlywedge x\neq\mathbf{0}$, whence, according to (\ref{T5-3}),
we conclude $y\in H$. This completes the proof of (\ref{53}).

To prove (\ref{56}) observe that for $x\in H$, by the
$\curlywedge$-quasi-stability, we also have $x[x\ldots
x]\curlywedge x\in H$, which in view of $x[x\ldots x]\curlywedge
x\leqslant x[x\ldots x]$ and (\ref{53}) implies $x[x\ldots x]\in
H$. So, (\ref{56}) is valid too.

Now, we shall verify (\ref{57}). Let $x\in H$. If $\mathbf{0}\in
H$ then, as it was proved in the proof of Theorem~\ref{T-5a},
$H=G$. Thus, in this case, $R_ix\in H$ for every
$i\in\overline{1,n}$. If $\mathbf{0}\not\in H$, then, evidently,
$x\neq\mathbf{0}$. Consequently, $R_ix\neq\mathbf{0}$ for every
$i\in\overline{1,n}$, because in the opposite case from
$\mathbf{A_2}$ it follows $x=\mathbf{0}$. This, together with
(\ref{56}), gives $x[x\ldots x]\neq\mathbf{0}$ and $R_ix[x\ldots
x]\neq\mathbf{0}$. Now, applying $\mathbf{A_6}$, we obtain
$\mathbf{0}\neq R_i((R_kx)[x\ldots x])$, whence, by (\ref{T5-1}),
we deduce $(R_kx)[x\ldots x]\neq\mathbf{0}$ for every
$k\in\overline{1,n}.$ Since $(R_kx)[x\ldots x]\leqslant x$, we
have
\[
(R_kx)[x\ldots x]\curlywedge x=(R_kx)[x\ldots x]\neq\mathbf{0},
\]
which, by (\ref{T5-3}), implies $R_kx\in H$ for all
$k\in\overline{1,n}$. The condition (\ref{57}) is proved.

The proof of (\ref{59}) is similar to the proof of (\ref{57}).
Namely, let $x[y\ldots y]\neq\mathbf{0}$ for some $x,y\in G$. If
$\mathbf{0}\in H$, then, as in the previous case, $H=G$. Hence
$R_ix\in H$ for every $i\in\overline{1,n}$. If $\mathbf{0}\not\in
H$, then $R_ix[y\ldots y]\neq\mathbf{0}$ for all
$i\in\overline{1,n}$, because in the case $R_ix[y\ldots y]=
\mathbf{0}$, by $\mathbf{A_2}$, we obtain $x[y\ldots y]=x[y\ldots
y][R_1x[y\ldots y]\ldots R_nx[y\ldots y]]=\mathbf{0}$ which
contradicts to our assumption. Next, applying $\mathbf{A_6}$, we
get $\mathbf{0}\neq R_ix[y\ldots y]=R_i((R_kx)[y\ldots y])$,
whence we deduce $(R_kx)[y\ldots y]\neq\mathbf{0}$ for each
$k\in\overline{1,n}$. In fact, from the above, for $(R_kx)[y\ldots
y]=\mathbf{0}$ it follows $R_i\mathbf{0}\neq\mathbf{0}$. This
contradicts to (\ref{T5-1}). Further, $(R_kx)[y\ldots y]\leqslant
y$ gives
\[
(R_kx)[y\ldots y]\curlywedge y=(R_kx)[y\ldots
y]\neq\mathbf{0},
\]
whence, by (\ref{T5-3}), we obtain $R_kx\in H$. This completes the
proof of (\ref{59}).

If $x[y\ldots y]=y\neq\mathbf{0}$ for some $x,y\in G$, then
$x[y\ldots y]\curlywedge y=x[y\ldots y]\neq\mathbf{0}$, whence,
according to (\ref{T5-3}), we have $x\in H$. This proves
(\ref{58}).

Now let $x\sqsubset y$ for some $x\in H$ and $y\in G$. Then,
obviously, $R_ix\leqslant R_iy$ for each $i\in\overline{1,n}$.
From this, applying (\ref{57}) and (\ref{53}), we obtain $R_iy\in
H$. So, (\ref{60}) is valid too.

At last, let $\mathbf{0}\not\in H$ and $x\sqsubset\mathbf{0}$.
Then $\mathbf{0}\leqslant x$, by (\ref{st-1a}), and
$\mathbf{0}=R_i\mathbf{0}$, by (\ref{T5-1}). Thus
$\mathbf{0}=R_i\mathbf{0}\leqslant R_ix$ for each
$i\in\overline{1,n}$. But from $x\sqsubset\mathbf{0}$ we have also
$R_ix\leqslant R_i\mathbf{0}=\mathbf{0}$. Therefore
$R_ix=\mathbf{0}$ for every $i\in\overline{1,n}$. Consequently,
$x=x[R_1x\ldots R_nx]=x[\mathbf{0}\ldots\mathbf{0}]=\mathbf{0}$.
This completes the proof of (\ref{61}) and the proof of
Theorem~\ref{T-prop}.
\end{proof}

\begin{minipage} {60mm}
\begin{flushleft}
Dudek~W. A. \\
 Institute of Mathematics and Computer Science \\
 Wroclaw University of Technology \\
 50-370 Wroclaw \\
 Poland \\
 E-mail: dudek@im.pwr.wroc.pl
\end{flushleft}
\end{minipage}
\hfill
\begin{minipage} {60mm}
\begin{flushleft}
 Trokhimenko~V. S. \\
 Department of Mathematics \\
 Pedagogical University \\
 21100 Vinnitsa \\
 Ukraine \\
 E-mail: vtrokhim@sovamua.com
 \end{flushleft}
 \end{minipage}

\end{document}